# LIMITING VELOCITY OF HIGH-DIMENSIONAL RANDOM WALK IN RANDOM ENVIRONMENT

By Noam Berger

*University of California, Los Angeles*

We show that random walk in uniformly elliptic i.i.d. environment in dimension $\geq 5$ has at most one non zero limiting velocity. In particular this proves a law of large numbers in the distributionally symmetric case and establishes connections between different conjectures.

**1. Introduction.** Let $d \geq 1$. A random walk in random environment (RWRE) on $\mathbb{Z}^d$ is defined as follows: Let $\mathcal{M}^d$ denote the space of all probability measures on $\{\pm e_i\}_{i=1}^d$ and let $\Omega = (\mathcal{M}^d)^{\mathbb{Z}^d}$. An *environment* is a point $\omega \in \Omega$. Let $P$ be a probability measure on $\Omega$. For the purposes of this paper, we assume that $P$ is an i.i.d. measure, that is,

$$P = Q^{\mathbb{Z}^d}$$

for some distribution $Q$ on $\mathcal{M}^d$ and that $P$ is *uniformly elliptic*, that is, there exists $\varepsilon > 0$ such that (s.t.) for every $e \in \{\pm e_i\}_{i=1}^d$,

$$Q(\{d : d(e) < \varepsilon\}) = 0.$$

For an environment $\omega \in \Omega$, the *random walk* on $\omega$ is a time-homogenous Markov chain with transition kernel

$$P_\omega(X_{n+1} = z + e | X_n = z) = \omega(z, e).$$

The *quenched law* $P_\omega^z$ is defined to be the law on $(\mathbb{Z}^d)^{\mathbb{N}}$ induced by the kernel $P_\omega$ and $P_\omega^z(X_0 = z) = 1$. We let $\mathbf{P} = P \otimes P_\omega^0$ be the joint law of the environment and the walk, and the *annealed* law is defined to be its marginal

$$\mathbb{P} = \int_\Omega P_\omega^0 \, dP(\omega).$$

We consider the limiting velocity

$$v = \lim_{n \to \infty} \frac{X_n}{n}.$$

Based on the work of Zerner [5] and Sznitman and Zerner [3], we know that $v$ exists $\mathbb{P}$-a.s. Furthermore, there is a set $A$ of size at most 2 such that almost surely $v \in A$.

Zerner and Merkl [6] proved that in dimension 2 a 0–1 law holds and therefore the set $A$ is of size 1, that is, a law of large numbers holds, in dimension 2 (see also [2] for a continuous version).

The main result of this paper is the following:

THEOREM 1.1. *For $d \geq 5$, there is at most one nonzero limiting velocity; that is, if $A = \{v_1, v_2\}$ with $v_1 \neq v_2$ and $v_1 \neq 0$, then $v_2 = 0$.*

Theorem 1.1 has the following immediate corollary:

COROLLARY 1.2. *For $d \geq 5$, if $Q$ is distributionally symmetric, then the limiting velocity is an almost sure constant.*

REMARK ABOUT CONSTANTS. As is common in most of the RWRE literature, the value of the constant $C$ may vary from line to line. In addition, $C$ may implicitly depend on variables that are kept constant throughout the entire calculation, in particular the dimension $d$ or the distribution $Q$.

**2. Backward path—Construction.** In this section we describe the backward path, the main object studied in this paper. The backward path is, roughly speaking, a path of the RWRE from $-\infty$ through the origin to $+\infty$. Below we define it. In Section 3 we prove some basic facts about it. Note that the backward path appears, though implicitly, in [1] and [4].

Throughout the paper we are assuming, for contradiction, that two different nonzero limiting velocities $v_1$ and $v_2$ exist. Assume without loss of generality that $\langle \ell, v_1 \rangle > 0$ for $\ell = e_1$. We let $A_\ell$ be the event that the walk is transient in the direction $\ell$, that is,

$$A_\ell = \left\{ \lim_{n \to \infty} \langle X_n, \ell \rangle = \infty \right\}.$$

By our assumptions, $Q$ is a distribution on $\mathcal{M}^d$ s.t. both $\mathbf{P}(A_\ell)$ and $\mathbf{P}(A_{-\ell})$ are positive.

We say that $t$ is a regeneration time in the direction $\ell$ if:

1. $\langle X_s, \ell \rangle < \langle X_t, \ell \rangle$ for every $s < t$, and
2. $\langle X_s, \ell \rangle > \langle X_t, \ell \rangle$ for every $s > t$.

REMARK. Note that in the special case of $\ell$ being a coordinate vector this simple definition coincides with the more complex definition of a regeneration time from [3].

For every $L > 0$, let $\mathcal{K}_L = \{z | 0 \leq \langle z, \ell \rangle < L\}$.

Let $t_1$ be the first regeneration time (if one exists), let $t_2$ be the second (if exists), and so on. If $t_{n+1}$ exists, let $L_n = \langle X_{t_{n+1}}, \ell \rangle - \langle X_{t_n}, \ell \rangle$, let

$$W_n : \mathcal{K}_{L_n} \to \mathcal{M}^d$$

be

$$W_n(z) = \omega(z + X_{t_n}),$$

let $u_n = t_{n+1} - t_n$ and let $K_n : [0, u_n] \to \mathbb{Z}^d$ be $K_n(t) = X_{t_n+t} - X_{t_n}$. We let $S_n$, the $n$th regeneration slab, be the ensemble $S_n = \{L_n, W_n, u_n, K_n\}$.

In [3] Sznitman and Zerner proved that on the event $A_\ell$, almost surely there are infinitely many regeneration times, and, furthermore, that the regeneration slabs $\{S_i\}_{i=1}^\infty$ form an i.i.d. process. Let $\lambda = \lambda_\ell$ be the distribution of $S_1$ conditioned on $A_\ell$.

We now construct an environment and a doubly infinite path in that environment. Let $\{S_n\}_{n \in \mathbb{Z}}$ be i.i.d. regeneration slabs sampled according to $\lambda$.

We now want to glue the regeneration slabs to each other. Let $Y_0 = 0$, and define, inductively, $Y_{n+1} = Y_n + K_n(u_n)$ for $n \geq 0$ and $Y_{n-1} = Y_n - K_{n-1}(u_{n-1})$ for $n \leq 0$. Almost surely $\mathbb{Z}^d$ is the disjoint union of the sets $Y_n + \mathcal{K}_{L_n}$. For every $z \in \mathbb{Z}^d$ let $n(z)$ be the unique $n$ such that $z \in Y_n + \mathcal{K}_{L_n}$. Let $\omega$ be the environment

$$\omega(z) = W_{n(z)}(z - Y_{n(z)}).$$

Let $\mathcal{T} \subseteq \mathbb{Z}^d$ be

$$\mathcal{T} = \bigcup_{n=-\infty}^{\infty} (Y_n + K_n[0, u_n]).$$

Let $\mu$ be the joint distribution of $\omega$ and $\mathcal{T}$. $\mathcal{T}$ is called the *backward path in direction* $\ell$. We let $\tilde{\mu}$ be the marginal distribution of $\omega$ in $\mu$.

**3. Backward path—Basic properties.** In this section we prove two simple properties of the measure $\mu$.

PROPOSITION 3.1. *There exists a coupling $\tilde{P}$ on $\Omega \times \Omega \times \{0,1\}^{\mathbb{Z}^d}$ with the distribution of $\omega, \tilde{\omega}, \mathcal{T}$ satisfying:*

1. *$\omega$ is distributed according to $P$.*
2. *$(\tilde{\omega}, \mathcal{T})$ is distributed according to $\mu$.*

3. $\tilde{P}$-almost surely, $\omega(z) = \tilde{\omega}(z)$ for every $z \in \mathbb{Z}^d \setminus \mathcal{T}$.
4. $\omega$ and $\mathcal{T}$ are independent.

PROPOSITION 3.2. *Let $\tilde{\omega}$ be an environment sampled according to $\tilde{\mu}$, and let $\{X_n\}$ be a random walk on that environment. Then almost surely $\{X_n\}$ is transient in the direction $\ell$.*

Both Proposition 3.1 and Proposition 3.2 follow from the fact that the $\tilde{\mu}$-environment around zero is similar to the $P$-environment around the location of the walker at a large regeneration time. More precisely, let $\omega, \{X_n\}$ be sampled according to $\mathbf{P}$ conditioned on the event $\forall_{n>0}(\langle X_n, \ell \rangle > 0) \cap A_\ell$, which is an event of positive probability. Let $t_1, t_2, \ldots$ be the regeneration times. (Note that we conditioned on transience in the $\ell$ direction, and therefore infinitely many regeneration times exist.) Let $\omega_i$ be the environment defined by $\omega_i(z) = \omega(z + X_{t_i})$ and let $\mathcal{T}_i \subseteq \mathbb{Z}^d$ be defined as $\mathcal{T}_i = \{X_t - X_{t_i} | t \geq 0\}$.

For $X \in \mathbb{Z}^d$ let $\mathcal{H}(X)$ be the half-space
$$\mathcal{H}(X) = \{z \mid \langle z, \ell \rangle \geq \langle X, \ell \rangle\}.$$

LEMMA 3.3. *For every $i$, the distribution of*
$$\{-X_{t_i}; \; \mathcal{T}_i \cap \mathcal{H}(-X_{t_i}); \; \omega_i|_{\mathcal{H}(-X_{t_i})}\} \tag{1}$$
*is the same as the distribution of*
$$\{Y_{-i}; \; \mathcal{T} \cap \mathcal{H}(Y_{-i}); \; \tilde{\omega}|_{\mathcal{H}(Y_{-i})}\}. \tag{2}$$

PROOF. Let $\tilde{\mathbf{P}}$ be $\mathbf{P}$ conditioned on the event $\forall_{n>0}(\langle X_n, \ell \rangle > 0) \cap A_\ell$. By Theorem 1.4 of [3], the distribution of
$$\{\omega|_{\mathcal{H}(0)}, \{X_t | t \geq 0\}\}$$
according to $\tilde{\mathbf{P}}$ is the same as the distribution of
$$\{\tilde{\omega}|_{\mathcal{H}(0)}, \mathcal{T} \cap \mathcal{H}(0)\}$$
according to $\mu$. The lemma now follows since the sequence $\{S_n\}_{n \in \mathbb{Z}}$ is i.i.d. □

We can now prove Propositions 3.1 and 3.2.

PROOF OF PROPOSITION 3.2. Let $B$ be the event that the walk is transient in the direction of $\ell$ and never exits the half-space $\mathcal{H}(0)$, that is,
$$B = A_\ell \cap \{\forall_t X_t \in \mathcal{H}(0)\}.$$

For a configuration $\omega$ and $z \in \mathbb{Z}^d$, let
$$R_\omega(z) = P_\omega^z(B).$$

Note that $R_\omega(z)$ depends only on $\omega|_{\mathcal{H}(0)}$, so by the Markov property
$$\mathbf{P}_\omega^{X_0}(B|X_1, X_2, \ldots, X_t) = R_\omega(X_t) \cdot \mathbf{1}_{X_1,\ldots,X_t \in \mathcal{H}(0)}.$$

In addition, $B \in \sigma(X_1, X_2, \ldots)$ and therefore almost surely
$$\lim_{t \to \infty} R_\omega(X_t) \geq \mathbf{1}_B.$$

In particular, $\tilde{\mathbf{P}}$-almost surely,
$$\lim_{t \to \infty} R_\omega(X_t) = 1,$$

and for the subsequence of regeneration times we get that $\tilde{\mathbf{P}}$-almost surely

(3) $$\lim_{n \to \infty} R_\omega(X_{t_n}) = 1,$$

and using the bounded convergence theorem, for
$$R_n = \mathbf{E}_{\tilde{\mathbf{P}}}(R_\omega(X_{t_n}))$$

we get

(4) $$\lim_{n \to \infty} R_n = 1.$$

Let $\{\tilde{\omega}, \mathcal{T}, \{Y_n\}\}$ be sampled according to $\mu$ and let $X_n$ be a random walk on the environment $\tilde{\omega}$, which is independent of $\{\mathcal{T}, \{Y_n\}\}$ conditioned on $\tilde{\omega}$. Let $B_N$ be the event
$$\lim_{n \to \infty} \langle X_n, \ell \rangle = \infty \quad \text{and} \quad \forall n \langle X_n, \ell \rangle \geq \langle Y_{-N}, \ell \rangle.$$

Then by Lemma 3.3

(5) $$(\mu \otimes P_{\tilde{\omega}}^0)(B_n) = R_n.$$

Remembering that
$$A_\ell = \bigcup_{n=1}^{\infty} B_n$$

we get from (5) that
$$(\mu \otimes P_{\tilde{\omega}}^0)(A_\ell) = \lim_{n \to \infty} R_n = 1,$$

as desired. $\square$

PROOF OF PROPOSITION 3.1. We define the coupling on every regeneration slab. Let $\tilde{\lambda}$ be the distribution on $\tilde{S} = \{L, W, \tilde{W}, u, K\}$ so that $\{L, \tilde{W}, u, K\}$ is distributed according to $\lambda$ and $W$ is defined as follows:

$$W(z) = \begin{cases} \tilde{W}(z), & \text{if } z \notin K([0, u]), \\ \psi(z), & \text{if } z \in K([0, u]), \end{cases}$$

where $\psi : \mathbb{Z}^d \to \mathcal{M}$ is sampled according to $P$, independently of $\{L, \tilde{W}, u, K\}$.

CLAIM 3.4.  *Conditioned on $L$, the environment $W$ is i.i.d. with marginal distribution $Q$, and independent of $u$ and $K$.*

We now sample the environments and the path as we did in Section 2: Let $\{\tilde{S}_n\}_{n=-\infty}^{\infty}$ be i.i.d. regeneration slabs sampled according to $\tilde{\lambda}$. Let $Y_0 = 0$ and define, inductively, $Y_{n+1} = Y_n + K_n(u_n)$ for $n \geq 0$ and $Y_{n-1} = Y_n - K_{n-1}(u_{n-1})$ for $n \leq 0$. Almost surely $\mathbb{Z}^d$ is the disjoint union of the sets $Y_n + \mathcal{K}_{L_n}$. For every $z \in \mathbb{Z}^d$ let $n(z)$ be the unique $n$ such that $z \in Y_n + \mathcal{K}_{L_n}$. We let $\omega$ be the environment

$$\omega(z) = W_{n(z)}(z - Y_{n(z)}),$$

we let $\tilde{\omega}$ be the environment

$$\tilde{\omega}(z) = \tilde{W}_{n(z)}(z - Y_{n(z)}),$$

and take $\mathcal{T} \subseteq \mathbb{Z}^d$ to be

$$\mathcal{T} = \bigcup_{n=-\infty}^{\infty} (Y_n + K_n[0, u_n]).$$

Clearly, $\{\tilde{\omega}, \mathcal{T}\}$ is distributed according to $\mu$ and $\omega$ and $\tilde{\omega}$ agree on $\mathbb{Z}^d - \mathcal{T}$. Therefore all we need to show is that $\omega$ is distributed according to $P$ and is independent of the path $\mathcal{T}$. This follows from Claim 3.4: conditioned on $\{u_n\}_{n=-\infty}^{\infty}$, $W$ is $P$-distributed and independent of the path $\mathcal{T}$. Therefore it is $P$-distributed and independent of the path $\mathcal{T}$ as we integrate over $\{u_n\}_{n=-\infty}^{\infty}$.  □

PROOF OF CLAIM 3.4. It is sufficient to show that conditioned on $L$, for every finite set $J = \{x_i : i = 1, \ldots, k\}$ with $J \subseteq \mathcal{K}_L$, the distribution of $\{W(x_i)\}_{x_i \in J}$ is i.i.d. with marginal $Q$ and independent of $u$ and $K$. This will follow if we prove that for every finite set $J = \{x_i | i = 1, \ldots, k\}$ with $J \subseteq \mathcal{K}_L$, conditioned on $L$, on $K$ and $u$ and on the event $J \cap K[0, u] = \varnothing$, the distribution of $\{\tilde{W}(x_i)\}_{x_i \in J}$ is i.i.d. with marginal $Q$.

To this end, fix $J$ and note that for every finite set $U$ that is disjoint of $J$, the event $\{K[0, u] = U\}$ is independent of $\{\tilde{W}(x_i)\}_{x_i \in J}$. Therefore, conditioned on the event $\{K[0, u] = U\}$ (and thus implicitly conditioning on $K$ and $u$), the distribution of $\{\tilde{W}(x_i)\}_{x_i \in J}$ is i.i.d. with marginal $Q$. By integrating with respect to $U$ we get that $\{W(x_i)\}_{x_i \in J}$ is $Q$-distributed, and by the fact that it was $Q$-distributed conditioned on $K$ and $u$ we get the independence. □

**4. Intersection of paths.** In this section we will see some interaction between the backward path and the path of an independent random walk.

Let $Q$ be a uniformly elliptic distribution so that $0 < \mathbf{P}(A_\ell) < 1$ and let $(\omega, \tilde{\omega}, \mathcal{T})$ be as in Proposition 3.1. Let $z_0$ be an arbitrary point in $\mathbb{Z}^d$, and let $\{X_i\}_{i=1}^\infty$ be a random walk on the configuration $\omega$ starting at $z_0$, such that:

1. $\{X_i\}$ is conditioned on the (positive probability) event that $\lim_{i \to \infty} \langle X_i, \ell \rangle = -\infty$.
2. Conditioned on $\omega$, $\{X_i\}_{i=1}^\infty$ is independent of $\tilde{\omega}$ and $\mathcal{T}$.

The purpose of this section is the following easy lemma:

LEMMA 4.1. *Under the conditions stated above, almost surely there exist infinitely many values of $i$ such that $X_i \in \mathcal{T}$.*

We will prove that almost surely there exists one such value of $i$. The proof that infinitely many exist is very similar but requires a little more care, and for the purpose of proving the main theorem of this paper one such $i$ is sufficient.

PROOF. We need to show that

(6) $$(\tilde{P} \otimes P_\omega^{z_0})\left(\lim_{i \to \infty} \langle X_i, \ell \rangle = -\infty \text{ and } \forall_i (X_i \notin \mathcal{T})\right) = 0.$$

In order to establish (6), let $\{Y_i\}_{i=1}^\infty$ be a random walk on the environment $\tilde{\omega}$, coupled to the rest of the probability space as follows:

Let
$$i_0 = \inf\{i : \omega(X_i) \neq \tilde{\omega}(X_i)\} \geq \inf\{i : X_i \in \mathcal{T}\}.$$
Now, for $i < i_0$, we define $Y_i = X_i$. For $i \geq i_0$, $Y_i$ is determined based on $Y_{i-1}$ according to $\tilde{\omega}(Y_{i-1})$ independently of $X_i$, $\omega$ and $\mathcal{T}$. Now, note that
$$\forall_i (X_i \notin \mathcal{T}) \implies i_0 = \infty \implies \forall_i (X_i = Y_i).$$
Therefore,
$$\left(\lim_{i \to \infty} \langle X_i, \ell \rangle = -\infty \text{ and } \forall_i (X_i \notin \mathcal{T})\right) \implies \lim_{i \to \infty} \langle Y_i, \ell \rangle = -\infty.$$
The proof is concluded if we remember that by Proposition 3.2,
$$(\tilde{P} \otimes P_{\tilde{\omega}}^{z_0})\left(\lim_{i \to \infty} \langle Y_i, \ell \rangle = -\infty\right) = 0.$$
□

## 5. Proof of main theorem.

LEMMA 5.1. *Let $d \geq 5$, and assume that the set $A$ of speeds contains two nonzero elements. Then there exists $z_0$ such that*

$$(\tilde{P} \otimes P_\omega^{z_0})\left(\lim_{i \to \infty} \langle X_i, \ell \rangle = -\infty \text{ and } \forall_i (X_i \notin \mathcal{T})\right) > 0.$$

PROOF. Let

$$\tilde{\mathcal{T}} = \{X_i : i = 1, 2, \ldots\}.$$

We use the following claim whose proof is deferred:

CLAIM 5.2. *Let $\tilde{B}$ be the event that $\langle X_i, \ell \rangle < \langle X_0, \ell \rangle$ for all $i > 0$. Note that $\tilde{B}$ has positive probability. Also, let $\mathcal{T}' = \mathcal{T} \cap \{z : \langle z, \ell \rangle \leq 0\}$. Then, if $A$ contains two distinct nonzero elements then*

(7) $$\sum_{z \in \mathbb{Z}^d} \tilde{P}(z \in \mathcal{T}')^2 < \infty$$

*and*

(8) $$\sum_{z \in \mathbb{Z}^d} \mathbb{P}^0(z \in \tilde{\mathcal{T}} | \tilde{B})^2 < \infty.$$

By Proposition 3.1, $\mathcal{T}'$ and $\tilde{\mathcal{T}}$ are independent random sets and therefore so are $\mathcal{T}'$ and $\tilde{\mathcal{T}} | \tilde{B}$. Therefore,

$$(\tilde{E} \otimes E_\omega^{z_0})(|\mathcal{T}' \cap \tilde{\mathcal{T}}||\tilde{B}) = \sum_{z \in \mathbb{Z}^d} \tilde{P}(z \in \mathcal{T}') \mathbb{P}^{z_0}(z \in \tilde{\mathcal{T}} | \tilde{B})$$

$$= \sum_{z \in \mathbb{Z}^d} \tilde{P}(z \in \mathcal{T}') \mathbb{P}^0(z - z_0 \in \tilde{\mathcal{T}} | \tilde{B}),$$

with the last equality following from translation invariance of the annealed measure. Let

$$M = \sum_{z \in \mathbb{Z}^d} \tilde{P}(z \in \mathcal{T}')^2$$

and

$$\tilde{M} = \sum_{z \in \mathbb{Z}^d} \mathbb{P}^0(z \in \tilde{\mathcal{T}} | \tilde{B})^2,$$

let $\lambda$ be so small that $\lambda M + \lambda \tilde{M} + \lambda^2 < 1$, and let $R$ be so large that

$$\sum_{\|z\| > R} \tilde{P}(z \in \mathcal{T}')^2 < \lambda \quad \text{and} \quad \sum_{\|z\| > R} \mathbb{P}^0(z \in \tilde{\mathcal{T}} | \tilde{B})^2 < \lambda.$$

Taking $z_0$ such that $\|z_0\| > 2R$ and $\langle z_0, \ell \rangle < 0$ we get, using Cauchy–Schwarz, that

$$(\tilde{E} \otimes E_\omega^{z_0})(|\mathcal{T}' \cap \tilde{\mathcal{T}}||\tilde{B})$$
$$= \sum_{z \in \mathbb{Z}^d} \tilde{P}(z \in \mathcal{T}')\mathbb{P}^0(z - z_0 \in \tilde{\mathcal{T}}|\tilde{B})$$
$$= \sum_{z \in B(0,R)} \tilde{P}(z \in \mathcal{T}')\mathbb{P}^0(z - z_0 \in \tilde{\mathcal{T}}|\tilde{B})$$
$$+ \sum_{z \in B(z_0,R)} \tilde{P}(z \in \mathcal{T}')\mathbb{P}^0(z - z_0 \in \tilde{\mathcal{T}}|\tilde{B})$$
$$+ \sum_{z \in \mathbb{Z}^d - B(0,R) - B(z_0,R)} \tilde{P}(z \in \mathcal{T}')\mathbb{P}^0(z - z_0 \in \tilde{\mathcal{T}}|\tilde{B})$$
$$\leq \lambda M + \lambda \tilde{M} + \lambda^2 < 1.$$

Therefore $\tilde{P} \otimes P_\omega^{z_0}(\mathcal{T}' \cap \tilde{\mathcal{T}} = \varnothing|\tilde{B}) > 0$. $P_\omega^{z_0}(\tilde{B}) > 0$ and by the choice of $z_0$, conditioned on $\tilde{B}$, $\mathcal{T}' \cap \tilde{\mathcal{T}} = \varnothing$ if and only if $\mathcal{T} \cap \tilde{\mathcal{T}} = \varnothing$. Therefore $\mathcal{T} \cap \tilde{\mathcal{T}}$ is empty with positive probability. $\square$

PROOF OF CLAIM 5.2. We will prove (7). Equation (8) follows from the exact same reasoning. First we get an upper bound on $\mu(Y_{-n} = z)$. The sequence $\{O_n = Y_{-n} - Y_{-n-1}\}$ is an i.i.d. sequence. Furthermore, due to ellipticity there exist $d$ linearly independent vectors $v_1, \ldots, v_d$ and $\varepsilon > 0$ such that for every $k = 1, \ldots, d$, and every $\delta \in \{+1, -1\}$,

$$\mu(O_1 = 2v_1 + \delta v_k) > \varepsilon.$$

($v_1$ is, approximately, in the direction of $\ell$, while the others are, approximately, orthogonal to $\ell$.)

Let

$$A = \{2v_1 + \delta v_k \mid k = 1, \ldots, d;\ \delta \in \{+1, -1\}\}$$

and let $p = \mu(O_1 \in A)$. Fix $n$, and let $E^{(n)}$ be the event that at least $\pi_n = \lceil \frac{1}{2} pn \rceil$ of the $O_i$'s, $i = 1, \ldots, n$, are in $A$. For every subset $H$ of $\{1, \ldots, n\}$ of size $\pi_n$, let $E_H^{(n)}$ be the event that the elements of $H$ are the smallest $\pi_n$ numbers $i$ such that $O_i \in A$. Then from heat kernel estimates for bounded i.i.d. random walks in $Z^d$ we get that for every $z \in \mathbb{Z}^d$,

$$\mu\left(\sum_{i \in H} O_i = z \Big| E_H^{(n)}\right) < Cn^{-d/2}.$$

Conditioned on $E_H^{(n)}$,

$$\sum_{i \in H} O_i \quad \text{and} \quad \sum_{i \notin H} O_i$$

are independent, so remembering that $Y_{-n} = \sum_{i=1}^{n} O_i$, we get that

$$\mu(Y_{-n} = z | E_H^{(n)}) < Cn^{-d/2}.$$

The events

$$\{E_H^{(n)} | H \subseteq [1, n]\}$$

are mutually exclusive and

$$\mu\left(\bigcup_H E_H^{(n)}\right) > 1 - e^{-Cn}.$$

Therefore, for every $n$ and $z \in \mathbb{Z}^d$,

(9) $$\mu(Y_{-n} = z) < Cn^{-d/2}.$$

Now, for every $n$ and $z \in \mathbb{Z}^d$, let $Q(z, n)$ be the probability that $z$ is visited during the $n$th regeneration, that is, between $Y_{1-n}$ and $Y_{-n}$. The $n$th regeneration is independent of $Y_{1-n}$, so

$$Q(z, n | Y_{1-n}) = Q(z - Y_{1-n}, 0).$$

The fact that the speed of the walk in direction $\ell$ is positive yields

(10) $$\sum_{z \in \mathbb{Z}^d} Q(z, 0) \leq E(\tau_2 - \tau_1) < \infty.$$

From (9) we get that

$$\sum_{z \in \mathbb{Z}^d} [\mu(Y_{-n} = z)]^2 \leq Cn^{-d/2}.$$

Combined with (10) and remembering that Young's inequality for convolution says that $\|f \star g\|_2 \leq \|f\|_2 \|g\|_1$ for all $f$ and $g$ (and noting that the next regeneration slab is independent of $Y_{1-n}$, and thus the result is a convolution), we get

$$\sum_{z \in \mathbb{Z}^d} [Q(z, n)]^2 \leq Cn^{-d/2}$$

or

(11) $$\sqrt{\sum_{z \in \mathbb{Z}^d} [Q(z, n)]^2} \leq Cn^{-d/4}.$$

Noting that

$$\mu(z \in \mathcal{T}') = \sum_{n=1}^{\infty} Q(z,n),$$

(11) and the triangle inequality tell us that

$$\sqrt{\sum_{z \in \mathbb{Z}^d} [\mu(z \in \mathcal{T}')]^2} \leq C \sum_{n=1}^{\infty} n^{-d/4}.$$

So for $d \geq 5$

$$\sum_{z \in \mathbb{Z}^d} [\mu(z \in \mathcal{T}')]^2 < \infty$$

as desired. $\square$

PROOF OF THEOREM 1.1. The theorem follows immediately from Lemma 4.1 and Lemma 5.1. $\square$

**Acknowledgments.** I thank G. Y. Amir, I. Benjamini, M. Biskup, N. Gantert, S. Sheffield and S. Starr for useful discussions. I thank O. Zeitouni for many important comments on a previous version of the paper. An anonymous referee is gratefully acknowledged for many important comments.

DEPARTMENT OF MATHEMATICS
UNIVERSITY OF CALIFORNIA, LOS ANGELES,
BOX 951555
LOS ANGELES, CALIFORNIA 90095-1555
USA
E-MAIL: berger@math.ucla.edu